\documentclass[12pt]{amsart}
\usepackage[T1]{fontenc}
\usepackage[utf8]{inputenc}
\usepackage{lmodern}
\usepackage[a4paper,margin=29mm]{geometry}
\usepackage{amsmath,amssymb,amsthm,mathtools}
\usepackage{microtype,enumitem}
\usepackage[hidelinks]{hyperref}
\hypersetup{
 pdftitle={On braided fusion categories of almost square-free dimension and pointed central de-equivariantizations},
 pdfauthor={Sebastian Burciu and Jingcheng Dong},
 pdfsubject={Fusion categories},
 pdfkeywords={braided fusion categories, almost square-free dimension, square-free dimension, equivariantization}}
\numberwithin{equation}{section}
\newtheorem{theorem}{Theorem}[section]
\newtheorem{proposition}[theorem]{Proposition}
\newtheorem{lemma}[theorem]{Lemma}
\newtheorem{corollary}[theorem]{Corollary}
\theoremstyle{definition}
\newtheorem{example}[theorem]{Example}
\theoremstyle{remark}
\newtheorem{remark}[theorem]{Remark}
\newcommand{\C}{\mathcal C}
\newcommand{\D}{\mathcal D}
\newcommand{\E}{\mathcal E}
\newcommand{\B}{\mathcal B}
\newcommand{\FPdim}{\operatorname{FPdim}}
\newcommand{\Rep}{\operatorname{Rep}}
\renewcommand{\Vec}{\operatorname{Vec}}
\newcommand{\sVec}{\operatorname{sVec}}
\newcommand{\Sem}{\operatorname{Sem}}
\newcommand{\Irr}{\operatorname{Irr}}
\newcommand{\Hom}{\operatorname{Hom}}
\newcommand{\id}{\operatorname{id}}
\newcommand{\one}{\mathbf1}
\newcommand{\Z}{\mathbb Z}
\newcommand{\Q}{\mathbb Q}
\newcommand{\br}{\mathrm{br}}
\DeclareMathOperator{\Fun}{Fun}
\DeclareMathOperator{\expn}{exp}

\setlist[enumerate]{label=\textup{(\roman*)},leftmargin=2.3em,itemsep=4pt}
\allowdisplaybreaks[1]
\title[Braided categories of almost square-free dimension]{On braided fusion categories of almost square-free dimension and pointed central de-equivariantizations}
\author{Sebastian Burciu}
\author{Jingcheng Dong}
\date{}
\subjclass[2020]{18M20, 16T05}
\keywords{Braided fusion category, almost square-free dimension, square-free dimension, equivariantization, fiber functor, quadratic form, M\"uger center}
\begin{document}

\begin{abstract}
Let \(p\) be an odd prime and let \(d\) be a square-free positive integer
coprime to \(p\). We determine the possible tensor structures of braided
fusion categories of Frobenius--Perron dimension \(p^n d\) for
\(0\le n\le4\). In odd dimension every such category is tensor equivalent
to the representation category of a finite group. In even dimension this
holds exactly when the M\"uger center has even dimension. In the remaining
case, where the category has even dimension and its center has odd
dimension, a semion splits off as a braided factor. Its complement is
tensor equivalent to the representation category of a group of odd order.
We reconstruct the group from a braided action on an odd-dimensional
pointed category and identify its extension class. Over \(\mathbb C\),
we construct a modular category of dimension \(p^5\) for each odd
prime \(p\) that shows the bound \(n\le4\) is sharp for tensor
reconstruction.
\end{abstract}
\maketitle
\section{Introduction}\label{asf:sec-introduction}

A positive integer is \emph{almost square-free} if it has the form
\(p^n d\), where \(p\) is prime, \(n\ge0\), and \(d\) is square-free
and coprime to \(p\). We ask when a braided fusion category of this dimension
is tensor equivalent to the representation category of a finite group.

Throughout, \(k\) is an algebraically closed field of characteristic zero,
and categories and functors are \(k\)-linear. Write \(\Vec\) for
finite-dimensional vector spaces over \(k\), \(\one\) for the tensor unit,
and \(\Rep(K)\) for finite-dimensional \(k\)-representations of a finite
group \(K\). A \emph{fiber functor} is a faithful exact strong tensor
functor to \(\Vec\), with no compatibility with a braiding required.
It realizes a fusion category as representations of a semisimple Hopf
algebra by reconstruction. A tensor equivalence with \(\Rep(K)\)
gives such a functor by composition with the usual forgetful functor.
We write \(\cong\) for isomorphism of groups or algebras,
\(\simeq_\otimes\) for tensor equivalence,
\(\simeq_{\br}\) for braided tensor equivalence, and \(\boxtimes\) for
the Deligne tensor product. In particular, a tensor equivalence with
\(\Rep(K)\) allows the transported braiding to differ from its usual
symmetric braiding.

We write \(\FPdim\) for Frobenius--Perron dimension. Throughout, the
dimension of a fusion category or an object means its Frobenius--Perron
dimension unless qualified. A fusion category
is \emph{integral} if the Frobenius--Perron dimension of each simple
object is an integer, and \emph{pointed} if all its simple objects are
invertible. For the braiding \(c\) of \(\C\), its \emph{M\"uger center}
\(\mathcal Z_2(\C)\) is the full subcategory of transparent objects,
namely those \(X\) such that \(c_{Y,X}c_{X,Y}=\id_{X\otimes Y}\)
for all \(Y\in\C\).

Dong and Dai proved that every braided fusion category of odd almost
square-free dimension is solvable. For \(n\le3\) they also realized
each such category as the representation category of a semisimple
quasitriangular Hopf algebra \cite{DD}. Here \emph{solvable} means
obtainable from \(\Vec\) by successive extensions and equivariantizations
by cyclic groups of prime order \cite[Proposition~4.4]{ENO}.
Theorem~\ref{asf:odd-asf-tensor} identifies the underlying tensor
category in odd dimension through \(n=4\). In even dimension it identifies
the obstruction to a realization as representations of a finite group.

Write \(C_m=\Z/m\Z\). For a finite group \(A\) and a normalized
\(3\)-cocycle \(\omega:A^3\to k^\times\), let \(\Vec_A^\omega\) denote
the category of \(A\)-graded vector spaces with associator \(\omega\)
on homogeneous degrees. We write \(\Vec_A\) for the trivial cocycle.
Fix \(i\in k\) with \(i^2=-1\). The semion categories \(\Sem_+\) and
\(\Sem_-\) are the pointed braided categories on \(C_2\) whose
nontrivial simple has self-braiding \(i\) and \(-i\), respectively.
\begin{theorem}
\label{asf:odd-asf-tensor}
Let \(p\) be an odd prime and let \(d\) be a square-free positive
integer coprime to \(p\). Let \(\C\) be a braided fusion category of
Frobenius--Perron dimension \(p^n d\), where \(0\le n\le4\).
Then \(\C\) is integral and satisfies the following conclusions.
\begin{enumerate}
\item If \(\FPdim(\C)\) is odd, then \(\C\simeq_\otimes\Rep(K)\)
      for a group \(K\) of order \(\FPdim(\C)\).
\item If \(\FPdim(\C)\) is even, exactly one of the following holds.
      \begin{enumerate}[label=\textup{(\alph*)}]
      \item The M\"uger center has even dimension and
            \(\C\simeq_\otimes\Rep(K)\) for a group \(K\)
            of order \(\FPdim(\C)\).
      \item The M\"uger center has odd dimension and
            \(\C\simeq_{\br}\Sem_\varepsilon\boxtimes\D\), where
            \(\varepsilon\in\{+,-\}\) and \(\D\) is a braided fusion
            category tensor equivalent to \(\Rep(H)\) for a group
            \(H\) of order \(\FPdim(\C)/2\). In particular,
            \begin{equation}\label{asf:eq-asf-semion-tensor}
             \C\simeq_\otimes
             \Vec_{C_2}^{\omega_{\mathrm{nt}}}\boxtimes\Rep(H),
             \qquad
             \omega_{\mathrm{nt}}(a,b,c)=(-1)^{abc}
             \quad (a,b,c\in\{0,1\}).
            \end{equation}
      \end{enumerate}
      In even dimension, the fiber functor criterion is
      \begin{equation}\label{asf:odd-asf-fiber-criterion}
       \C\text{ admits a fiber functor}
       \quad\Longleftrightarrow\quad
       2\mid\FPdim\bigl(\mathcal Z_2(\C)\bigr).
      \end{equation}
\end{enumerate}
\end{theorem}

A fusion category is \emph{group-theoretical} if it is categorically
Morita equivalent to a pointed fusion category. For \(0\le n\le4\),
Proposition~\ref{asf:central-pointed-reduction} gives group-theoreticality
under the hypotheses of Theorem~\ref{asf:odd-asf-tensor}. Solvability
in even dimension holds for every \(n\ge0\). Indeed, if \(d\) is even,
write \(p^n d=p^n 2(d/2)\). For \(n\ge1\) Dong's dimension theorem
\cite{D26} gives weak group-theoreticality. For \(n=0\) this also follows
from Proposition~\ref{asf:central-pointed-reduction}, since a pointed
category is nilpotent. In either case \(p^n d\equiv2\pmod4\), so
\cite[Lemma~4.3]{CCP} gives solvability.

Natale showed that every quasitriangular Hopf algebra of odd square-free
dimension is a group algebra \cite[Theorem~1.2]{Nat06}. Brugui\`eres
and Natale used this result to identify the underlying tensor category
in odd square-free dimension \cite[Theorem~1.1]{BN}.
Taking \(n=0\) in Theorem~\ref{asf:odd-asf-tensor}, with \(d=N\) square-free
and an odd prime \(p\nmid N\), also gives the even-dimensional case.

\begin{corollary}
\label{asf:squarefree-tensor}
Let \(\C\) be a braided fusion category of square-free dimension \(N\).
If \(N\) is odd, or if \(N\) and \(\FPdim(\mathcal Z_2(\C))\) are both
even, then \(\C\simeq_\otimes\Rep(K)\) for a group \(K\) of order \(N\).
If \(N\) is even and \(\FPdim(\mathcal Z_2(\C))\) is odd, then
\(\C\simeq_{\br}\Sem_\varepsilon\boxtimes\D\) for some
\(\varepsilon\in\{+,-\}\) and a braided fusion category \(\D\)
tensor equivalent to \(\Rep(H)\) for a group \(H\) of order \(N/2\). In this last case \(\C\) has no fiber
functor.
\end{corollary}

Bao, Yu, and Zheng prove that a braided finite tensor category over
\(\mathbb C\) of square-free dimension admitting a fiber functor is
semisimple and tensor equivalent to \(\Rep(K)\) for a finite group \(K\)
\cite[Corollary~3.8]{BYZ}. In the fusion case,
Corollary~\ref{asf:squarefree-tensor} gives the same conclusion, since
Lemma~\ref{asf:semion-associator} rules out a fiber functor on the
semion factor.
Remark~\ref{byz} explains how this also recovers the semisimple case
of their Hopf algebra theorem \cite[Theorem~1.1]{BYZ}.

Section~\ref{asf:sec-preliminaries} develops the central
de-equivariantization used in the proof. In
Section~\ref{asf:sec-reconstruction} we reconstruct a finite group from
a braided action on an odd-dimensional pointed category. The coherence
of the action determines the extension class of the group. We then show
that in any braided fusion category of dimension twice an odd integer,
a semion splits off exactly when the M\"uger center has odd dimension
(Theorem~\ref{asf:general-semion-splitting}).
Section~\ref{asf:subsec-main-proof} proves
Theorem~\ref{asf:odd-asf-tensor}. In
Section~\ref{asf:sharpness-exponent-five}, Theorem~\ref{asf:sharpness-theorem}
constructs, for each odd prime \(p\), a counterexample in dimension
\(p^5\). Taking its product with \(\Rep(C_2)\) or with a semion gives
counterexamples to the corresponding group-representation conclusions
in even dimension. The remaining
subsections illustrate the role of the coherence class and the prime
two, give applications to Hopf algebras, and describe the possible
nonpointed primary factor at exponent five
(Proposition~\ref{asf:exponent-five-obstruction}).

\section{Preliminaries and central de-equivariantization}
\label{asf:sec-preliminaries}

\subsection{Dimension and centralizers}
\label{asf:subsec-dimensions}

We follow \cite{EGNO}. Write \(\Irr(\C)\) for the isomorphism classes
of simple objects and \(\C_{\mathrm{pt}}\) for the fusion subcategory
generated by the invertible objects. For a simple \(X\), its fusion
matrix records the multiplicities of the simples in \(X\otimes Y\)
as \(Y\) runs through \(\Irr(\C)\). Its spectral radius is
\(\FPdim(X)\). The dimension of the category is
\begin{equation}\label{asf:eq-category-dimension}
 \FPdim(\C)=\sum_{X\in\Irr(\C)}\FPdim(X)^2.
\end{equation}
A fusion category is \emph{weakly integral} if \(\FPdim(\C)\) is an integer.

For a fusion subcategory \(\D\) of a braided fusion category \(\C\), its M\"uger centralizer
\(\D'_\C\) consists of the objects \(X\in\C\) satisfying
\(c_{Y,X}c_{X,Y}=\id_{X\otimes Y}\) for every \(Y\in\D\). Thus
\(\C'_\C=\mathcal Z_2(\C)\). We call \(\C\) \emph{nondegenerate} when
\(\mathcal Z_2(\C)\simeq_{\br}\Vec\), and \emph{slightly degenerate}
when \(\mathcal Z_2(\C)\simeq_{\br}\sVec\). Here \(\sVec\) is the
category of super vector spaces with the sign braiding. Its odd line
is the \emph{fermion}, an invertible simple of order two and
self-braiding \(-1\). We use the centralizer dimension formula
\cite[Theorem~3.14, Eq.~(37)]{DGNO}
\begin{equation}\label{asf:eq-centralizer-dimension}
 \FPdim(\D)\FPdim(\D'_\C)
 =\FPdim(\C)\FPdim\bigl(\D\cap\mathcal Z_2(\C)\bigr).
\end{equation}
In an integral fusion category, every fusion subcategory satisfies
\begin{equation}\label{asf:eq-subcategory-divisibility}
 \FPdim(\D)\mid\FPdim(\C),
\end{equation}
since the quotient is a rational algebraic integer
\cite[Theorem~7.17.6]{EGNO}. A pivotal structure identifies the
identity functor monoidally with the double-dual functor. It is
spherical when the left and right traces agree. We equip each
integral braided category with its canonical positive spherical
structure, whose categorical dimensions agree with its Frobenius--Perron
dimensions \cite[Corollary~2.24 and Remark~2.25]{DGNO}. A nondegenerate
braided category with this structure is modular, with the ribbon twist
induced by its braiding and pivotal structure. For a modular category,
the Frobenius--Schur exponent \(\operatorname{FSexp}(\C)\) is the order
of the twist matrix \(T=\operatorname{diag}(\theta_X)_{X\in\Irr(\C)}\),
where \(\theta_X\) is the ribbon twist on \(X\) and \(\theta_{\one}=1\)
\cite[Theorem~7.7]{NS}.

The adjoint subcategory \(\C_{\mathrm{ad}}\) is generated by the simple
summands of \(X\otimes X^*\) for all simple \(X\). A category is
\emph{nilpotent} if iterating this construction eventually gives \(\Vec\).
A fusion category is \emph{weakly group-theoretical} if it is
categorically Morita equivalent to a nilpotent fusion category
\cite[Definition~1.1]{ENO}. For a pointed category the adjoint
subcategory is \(\Vec\), so every group-theoretical category is weakly
group-theoretical.

The universal grading of \(\C\) has grading group \(U(\C)\) and
neutral component \(\C_{\mathrm{ad}}\). Every faithful group grading is
obtained from it via a quotient of \(U(\C)\). Its components have equal
Frobenius--Perron dimension, and for a modular category
\(|U(\C)|=\FPdim(\C_{\mathrm{pt}})\) \cite[\S2.1]{CCP}.

\begin{proposition}[{\cite[Corollary~3.2]{DN}}]
\label{asf:integrality-input}
A weakly integral braided fusion category \(\C\) with
\(4\nmid\FPdim(\C)\) is integral.
\end{proposition}
This applies to square-free dimension and to \(p^n d\) when
\(p\) is odd, \(d\) is square-free, and \(\gcd(p,d)=1\).

\subsubsection*{Numerical invariants over the ground field}
We will apply the complex modular classification of \cite{CPS} to
integral nondegenerate braided categories over \(k\). To justify this,
take a split skeletal presentation with one label for each simple object
and choose bases of its fusion spaces. The fusion multiplicities and the
matrices for the associator, unit constraints, braiding, rigidity and
positive spherical structure determine this presentation. Their entries,
together with the entries of the inverse isomorphisms, generate a finitely
generated subfield \(F\subset k\) over \(\Q\). Using the same labels and
multiplicity spaces over \(F\) gives a split semisimple model with
endomorphism ring \(F\) for each simple label. The coherence and
duality identities hold over \(F\) because they hold after the
inclusion \(F\hookrightarrow k\).

Since \(F\) is finitely generated over \(\Q\), it embeds in
\(\mathbb C\). Extending the model
along this embedding preserves the fusion matrices, hence the simple
Frobenius--Perron dimensions and pointedness. The spherical dimensions
are the same positive integers. The entries of the matrix of traces of
double braidings (the \(S\)-matrix) and the twist scalars are carried to
their images. Thus the determinant of \(S\) remains nonzero and the
order of each twist scalar is preserved. The complex model is modular,
and the pointed subcategory dimension and Frobenius--Schur exponent
computed there also hold over \(k\).

\subsection{Categorical actions and central de-equivariantization}
\label{asf:subsec-central}

A braided fusion subcategory is \emph{Tannakian} if it is braided
equivalent to \(\Rep(G)\) with its usual symmetric braiding for a finite
group \(G\). It is \emph{central} when it lies in the M\"uger center.

A categorical action of a finite group \(G\) on a fusion category \(\B\)
consists of tensor autoequivalences \(F_g\) and monoidal natural
isomorphisms \(\mu_{g,h}:F_gF_h\to F_{gh}\). We normalize the unit data
by \(F_e=\id\) and \(\mu_{e,g}=\mu_{g,e}=\id\). The action pentagon is
\begin{equation}\label{asf:eq-action-pentagon}
 \mu_{gh,l}\circ(\mu_{g,h}F_l)
 =\mu_{g,hl}\circ F_g(\mu_{h,l}).
\end{equation}
An equivariant object is an object \(V\) with isomorphisms
\(u_g:F_g(V)\to V\) satisfying
\begin{equation}\label{asf:eq-equivariance}
 u_g\circ F_g(u_h)=u_{gh}\circ\mu_{g,h,V},\qquad u_e=\id_V.
\end{equation}
These objects, with morphisms commuting with the equivariance maps,
form the equivariantization \(\B^G\). For a braided action, each
\(F_g\) is braided, the braiding descends to \(\B^G\), and objects
whose underlying objects are direct sums of \(\one\) form a central
Tannakian subcategory \(\Rep(G)\). Isomorphisms of two actions on the
same category are families of monoidal natural isomorphisms compatible
with the unit and the maps \(\mu_{g,h}\). More generally, an equivariant
tensor equivalence induces a tensor equivalence of equivariantizations,
which is braided when the given equivalence is braided
\cite[Sections~2.7 and~4.15]{EGNO}.

Every symmetric fusion category is of the form \(\Rep(L,u)\) for a
finite group \(L\) and a central element \(u\) with \(u^2=1\). Its
braiding is the sign braiding for the parity defined by \(u\), and its
unique maximal Tannakian subcategory is \(\Rep(L/\langle u\rangle)\)
\cite[Theorem~2.43 and proof of Corollary~2.50]{DGNO}. In particular, the
M\"uger center of \(\C\) has a unique maximal Tannakian subcategory.
We choose a symmetric equivalence and write
\begin{equation}\label{asf:eq-central-tannakian}
 \E_\C=\Rep(G)\subseteq\mathcal Z_2(\C).
\end{equation}

For any central Tannakian subcategory \(\Rep(G)\subseteq\mathcal Z_2(\C)\),
the de-equivariantization \(\C_G\) is the category of right modules in
\(\C\) over the regular algebra \(R=\Fun(G,k)\), with pointwise
multiplication and left translation. Its tensor product is the relative
tensor product over \(R\). The braiding descends to \(\C_G\), and
right translations induce a braided categorical action of \(G\).
The inverse equivalence and dimension formula
\cite[Proposition~4.22(i)--(ii) and Proposition~4.26]{DGNO}
give
\begin{equation}\label{asf:eq-central-equivariantization}
 \C\simeq_{\br}(\C_G)^G,\qquad
 \FPdim(\C_G)=\frac{\FPdim(\C)}{|G|}.
\end{equation}
The de-equivariantization \(\C_G\) is integral if and only if \(\C\) is integral
\cite[Corollary~4.27]{DGNO}. Applying the centralizer
correspondence \cite[Proposition~4.30(iii)]{DGNO} to
\(\C\) itself gives
\begin{equation}\label{asf:eq-center-quotient}
 \mathcal Z_2(\C_G)\simeq_{\br}\bigl(\mathcal Z_2(\C)\bigr)_G.
\end{equation}
Take \(\Rep(G)=\E_\C\) as in~\eqref{asf:eq-central-tannakian}.
If \(\mathcal Z_2(\C)\) is Tannakian, then \(\E_\C=\mathcal Z_2(\C)\),
so~\eqref{asf:eq-center-quotient} gives
\(\mathcal Z_2(\C_G)\simeq_{\br}\Vec\).
Otherwise \(\E_\C\) has index two in \(\mathcal Z_2(\C)\), and
\(\bigl(\mathcal Z_2(\C)\bigr)_G\) is a symmetric
non-Tannakian category of Frobenius--Perron dimension two.
Thus \(\mathcal Z_2(\C_G)\simeq_{\br}\sVec\)
\cite[Corollaries~2.50(ii) and~4.31(ii)]{DGNO}.

\subsection{Pointed categories, coefficient actions, and semions}
\label{asf:subsec-pointed-categories}

Every pointed fusion category over \(k\) is tensor equivalent to
\(\Vec_A^\omega\). Here \(A\) is the group of isomorphism classes of
simple objects under tensor product. Choosing representatives and
identifying their tensor products gives a scalar associator
\(\omega(a,b,c)\). The pentagon is the \(3\)-cocycle identity, and
changing these identifications changes \(\omega\) by a coboundary.
Thus \(\Vec_A^\omega\) and \(\Vec_{A'}^{\omega'}\) are tensor equivalent
exactly when some group isomorphism \(f:A\to A'\) satisfies
\(f^*[\omega']=[\omega]\) in \(H^3(A,k^\times)\), with trivial
coefficient action \cite[Proposition~2.6.1 and Remark~2.6.2]{EGNO}.

A braiding makes \(A\) abelian. If \(X_a\) represents \(a\in A\),
write \(q(a)\) for the scalar by which \(c_{X_a,X_a}\) acts on
\(X_a\otimes X_a\). The classification of pointed braided fusion
categories identifies their braided equivalence classes with isometry
classes of \emph{pre-metric groups} \((A,q)\)
\cite[Proposition~2.41]{DGNO}, \cite[Theorem~8.4.12]{EGNO}.
Here \(A\) is finite abelian, \(q:A\to k^\times\) satisfies
\(q(0)=1\) and \(q(-a)=q(a)\), and its polarization \(b_q\) below is
multiplicative in each variable. An isometry is a group isomorphism
preserving \(q\). We denote the corresponding category by \(\C(A,q)\).
Its double braiding on degrees \(a,b\) is
\begin{equation}\label{asf:eq-bicharacter}
 b_q(a,b)=q(a+b)q(a)^{-1}q(b)^{-1}.
\end{equation}
For a subgroup \(L\le A\), write
\(L^\perp=\{a\in A:b_q(a,x)=1\text{ for all }x\in L\}\).
In particular, \(A^\perp\) is the radical of \(b_q\). The M\"uger
center is the pointed subcategory on \(A^\perp\) with the restricted form.
We use the identities
\begin{equation}\label{asf:eq-quadratic-identities}
 b_q(a,a)=q(a)^2,\qquad q(ra)=q(a)^{r^2}\quad(r\in\Z).
\end{equation}
Distinct primary components of \(A\) are orthogonal since a mixed
bicharacter value has order dividing two coprime integers.
Orthogonal sums of pre-metric groups correspond to braided Deligne
products.
The term \emph{pre-metric} allows a nontrivial radical. Write
\(\expn(A)\) for the least common multiple of the orders of elements of
\(A\), \(O(A,q)\) for the group of isometries, and
\(\widehat A=\Hom(A,k^\times)\). For an action
\(\rho:G\to O(A,q)\), we use the induced \emph{left} action
\begin{equation}\label{asf:eq-coefficient-action}
 (g\cdot\chi)(a)=\chi(\rho_g^{-1}a)
 \qquad(\chi\in\widehat A,\ a\in A).
\end{equation}
For an abelian \(G\)-module \(M\) written multiplicatively,
\(Z^2(G,M)\) denotes the group of normalized \(2\)-cocycles,
that is, maps \(\eta:G\times G\to M\) satisfying
\begin{equation}\label{asf:eq-two-cocycle-definition}
 \eta(e,g)=\eta(g,e)=1,\qquad
 \eta(g,h)\eta(gh,l)=(g\cdot\eta(h,l))\eta(g,hl).
\end{equation}
The coboundary of a normalized cochain \(\lambda:G\to M\) is
\((\delta\lambda)(g,h)=\lambda(g)(g\cdot\lambda(h))\lambda(gh)^{-1}\).
Quotienting \(Z^2(G,M)\) by these coboundaries gives \(H^2(G,M)\).
For \(M=\widehat A\) with action~\eqref{asf:eq-coefficient-action},
it classifies group extensions with fixed kernel, quotient, and
conjugation action \cite[Theorem~6.6.3]{W}.

\begin{lemma}
\label{asf:semion-associator}
The semion categories
\(\Sem_+=\C(C_2,q_+)\) and \(\Sem_-=\C(C_2,q_-)\), with
\(q_+(1)=i\) and \(q_-(1)=-i\), are nondegenerate. Their common
underlying tensor category is \(\Vec_{C_2}^{\omega_{\mathrm{nt}}}\)
with associator
\begin{equation}\label{asf:eq-semion-cocycle}
 \omega_{\mathrm{nt}}(a,b,c)=(-1)^{abc}\qquad(a,b,c\in\{0,1\}).
\end{equation}
The class of this associator is nontrivial, and neither semion category
admits a fiber functor.
\end{lemma}
\begin{proof}
For \(G=C_2\) and coefficients in \(k^\times\), Braunling's formula
\cite[Theorem~7.1]{Br} gives the braiding
\(c_\varepsilon(a,b)=q_\varepsilon(1)^{ab}\) and the associator
\(\omega(a,b,c)=q_\varepsilon(1)^{2abc}=(-1)^{abc}\) for either sign.
See also Rowell--Stong--Wang \cite[\S5.3.1]{RSW} for the
explicit braiding and associator of \(\Sem_+\).
The double braiding of the nontrivial simple object with itself is
\(q_\varepsilon(1)^2=-1\), so both categories are nondegenerate.
Every normalized \(3\)-coboundary on \(C_2\) takes the value \(1\) at
\((1,1,1)\), whereas \(\omega_{\mathrm{nt}}(1,1,1)=-1\).
Thus \([\omega_{\mathrm{nt}}]\ne 0\), and the common underlying tensor
category admits no fiber functor by \cite[Example~5.1.3]{EGNO}.
\end{proof}

For comparison, \(\sVec\) has \(q(1)=-1\) and trivial associator.
For any finite abelian group \(A\), evaluation gives a tensor equivalence
\(\Vec_A\simeq_\otimes\Rep(\widehat A)\).

\section{Reconstruction from a pointed central de-equivariantization}
\label{asf:sec-reconstruction}

We recall two cohomological facts used in the reconstruction. The first
compares tensor constraints on braided functors. The second controls
the coherence isomorphisms of the group action in the coprime case.

\begin{lemma}
\label{asf:symmetric-cocycle}
Let \(A\) be a finite abelian group, and let
\(J:A\times A\to k^\times\) be a normalized \(2\)-cocycle for the
trivial action. If \(J(a,b)=J(b,a)\) for all \(a,b\in A\), then
\(J\) is a coboundary.
\end{lemma}
\begin{proof}
This is \cite[Exercise~8.4.13]{EGNO}. For the cohomological
description, see also Galindo--Morales \cite[\S2, (2-5)]{GM}.
Over our algebraically closed field, the map sending
\([J]\in H^2(A,k^\times)\) to the alternating bicharacter
\((a,b)\mapsto J(a,b)/J(b,a)\) is an isomorphism.
For symmetric \(J\) this bicharacter is trivial,
so \(J\) is a coboundary.
\end{proof}

\begin{lemma}
\label{asf:coprime-cohomology}
Let \(G\) be a finite group acting on a finite abelian group \(M\).
If \(\gcd(|G|,|M|)=1\), then \(H^2(G,M)=0\).
\end{lemma}
\begin{proof}
Since \(\gcd(|G|,|M|)=1\), multiplication by \(|G|\) is an
automorphism of \(M\). Hence \(H^2(G,M)=0\) by
\cite[Proposition~6.1.10]{W}.
\end{proof}

\subsection{An invariant model in odd order}

\begin{lemma}
\label{asf:odd-canonical-model}
Let \((A,q)\) be a pre-metric group with \(|A|\)
odd. For any integer
\(r\) with \(2r\equiv1\pmod{\expn(A)}\), the category \(\C(A,q)\) is
braided equivalent to \(\Vec_A\) with trivial associator and braiding
\begin{equation}\label{asf:eq-canonical-braiding}
 c_q(a,b)=b_q(a,b)^r\qquad(a,b\in A).
\end{equation}
This symmetric bicharacter is independent of \(r\), its double braiding
is \(b_q\), and every isometry of \((A,q)\) acts by strict braided regrading.
\end{lemma}
\begin{proof}
Put \(N=\expn(A)\). Every value of \(b_q\) is an \(N\)-th root of
unity. Thus \(c_q=b_q^r\) is a symmetric bicharacter independent of the
choice of inverse \(r\) of \(2\) modulo \(N\), and its double braiding
is \(b_q\). To compute its diagonal, fix \(a\in A\) of order \(m\).
From~\eqref{asf:eq-quadratic-identities},
\(q(a)^{m^2}=q(ma)=1\) and
\(q(a)^{2m}=b_q(a,ma)=1\). Since \(m\) is odd,
\(\gcd(m^2,2m)=m\). It follows that \(q(a)^m=1\), hence
\[
 c_q(a,a)=q(a)^{2r}=q(a).
\]
The classification by pre-metric groups
\cite[Theorem~8.4.12]{EGNO} now identifies \((\Vec_A,c_q)\) with
\(\C(A,q)\) by a braided equivalence inducing the identity on the
simple labels \(A\). Every isometry preserves \(b_q\) and \(c_q\),
so its regrading functor has identity tensor and unit constraints.
This is the invariant bicharacter construction of
\cite[\S5.1.1, Corollaries~1--2]{CGPW}, here written over \(k\).
\end{proof}

\subsection{Reconstructing the group and its representations}

\begin{theorem}
\label{asf:odd-pointed-equivariantization}
Let \((A,q)\) be a pre-metric group with \(|A|\) odd,
and let a finite group \(G\) act on
\(\mathcal P=\C(A,q)\) by braided tensor autoequivalences.
There is a finite group \(K\) such that
\begin{equation}\label{asf:eq-reconstructed-category}
 \mathcal P^G\simeq_\otimes\Rep(K),
 \qquad |K|=|A||G|.
\end{equation}
More precisely, the action induces a homomorphism
\(\rho:G\to O(A,q)\) and an extension
\begin{equation}\label{asf:eq-reconstruction-extension}
 1\longrightarrow\widehat A\longrightarrow K\longrightarrow G
 \longrightarrow1,
 \qquad
 (g\cdot\chi)(a)=\chi(\rho_g^{-1}a).
\end{equation}
For fixed \(\rho\), the isomorphism classes of braided categorical actions
of \(G\) on \(\mathcal P\) inducing \(\rho\) are parametrized by
\(H^2(G,\widehat A)\). An isomorphism of two such actions is a family of
monoidal natural isomorphisms \(F_g\Rightarrow F'_g\) compatible with
the unit and composition constraints of the two actions.
The parameter is the extension class in
\eqref{asf:eq-reconstruction-extension}.
In particular, if \(\gcd(|G|,|A|)=1\), then one may take
\(K=\widehat A\rtimes G\).
\end{theorem}
\begin{proof}
A braided autoequivalence permutes the invertible simple objects,
preserves their tensor products, and preserves their scalar
self-braidings. It therefore induces an isometry of \((A,q)\).
The coherence isomorphisms imply
\(\rho_g\rho_h=\rho_{gh}\), so \(\rho\) is a group homomorphism.

Use the model in Lemma~\ref{asf:odd-canonical-model}.
Fix \(g\in G\). Since the category is semisimple and the simple
endomorphism spaces are \(k\), an autoequivalence inducing \(\rho_g\)
is naturally isomorphic as a linear functor to regrading by \(\rho_g\).
Choose this isomorphism on the tensor unit to identify the unit constraint
with the identity, and transport the tensor structure along it.
On simple objects \(a,b\), the tensor constraint is a scalar
\(J_g(a,b)\). The monoidal functor identity says that \(J_g\) is
a normalized \(2\)-cocycle on \(A\).
The braided functor identity reads
\begin{equation}\label{asf:eq-braided-tensor-constraint}
 c_q(\rho_g a,\rho_g b)J_g(b,a)
 =
 J_g(a,b)c_q(a,b).
\end{equation}
Because \(\rho_g\) preserves \(c_q\), this implies
\(J_g(a,b)=J_g(b,a)\).
The symmetric cocycle lemma (Lemma~\ref{asf:symmetric-cocycle})
shows that the tensor structure can be removed by a monoidal
natural isomorphism.

Transport the action along these isomorphisms. For each \(g\in G\),
the functor \(F_g\) is now the strict braided tensor autoequivalence
of \(\C(A,q)\) realizing the isometry \(\rho_g\) in the model of
Lemma~\ref{asf:odd-canonical-model}. It sends the simple object of
degree \(a\) to the simple object of degree \(\rho_g(a)\). Thus for
an \(A\)-graded vector space \(V\) we have
\begin{equation}\label{asf:eq-strict-action}
 (F_gV)_a=V_{\rho_g^{-1}a}.
\end{equation}
The composites \(F_gF_h\) agree with \(F_{gh}\), but the action coherence
maps need not be identities. A monoidal natural automorphism of a
regrading functor is a scalar on each simple output degree.
Monoidality means that these scalars multiply when output degrees
are added. They therefore form a character of \(A\).
Write the component of \(\mu_{g,h}\) in output degree \(a\) as
\(\eta(g,h)(a)\), where \(\eta(g,h)\in\widehat A\).

Substitution into the action
pentagon~\eqref{asf:eq-action-pentagon} gives
\begin{equation}\label{asf:eq-action-two-cocycle}
 \eta(g,h)\eta(gh,l)
 =
 (g\cdot\eta(h,l))\eta(g,hl).
\end{equation}
To see the action on the right, the component of
\(F_g(\mu_{h,l})\) in output degree \(a\) is
\(\eta(h,l)(\rho_g^{-1}a)\).
This uses the coefficient action~\eqref{asf:eq-coefficient-action}.
The unit normalization gives \(\eta(e,g)=\eta(g,e)=1\). Conversely,
any normalized cocycle defines coherence maps on these regrading
functors. More explicitly, for two such actions with cocycles
\(\eta,\eta'\), an action isomorphism from the first to the second
has components
\(\lambda(g)\in\widehat A\) on output degrees. Compatibility with
composition reads
\[
 \lambda(gh)\eta(g,h)
   =\eta'(g,h)\lambda(g)(g\cdot\lambda(h)).
\]
Hence \(\eta'=\eta(\delta\lambda)^{-1}\), and any such normalized
cochain gives an isomorphism. This proves the stated parametrization
for the fixed action \(\rho\) on the labelled group \(A\).

The extension with class \([\eta]\) in
\eqref{asf:eq-reconstruction-extension} has underlying set
\(\widehat A\times G\) and multiplication
\cite[Theorems~6.6.3 and~6.6.7]{W}
\begin{equation}\label{asf:eq-extension-multiplication}
 (\chi,g)(\psi,h)
 =
 \bigl(\chi(g\cdot\psi)\eta(g,h),gh\bigr).
\end{equation}
In particular, \(|K|=|\widehat A||G|=|A||G|\).

It remains to construct the tensor equivalence.
Let \((V,(u_g))\) be an equivariant object and regard
\(u_g:F_g(V)\to V\) as a linear operator \(T_g\) on the underlying
vector space. By~\eqref{asf:eq-strict-action},
\(T_g(V_a)\subseteq V_{\rho_g a}\).
For \(\chi\in\widehat A\), define \(D_\chi\) by
\(D_\chi|_{V_a}=\chi(a)\id_{V_a}\).
The grading and the equivariance
identity~\eqref{asf:eq-equivariance} give
\begin{equation}\label{asf:eq-representation-operators}
 T_gD_\chi T_g^{-1}=D_{g\cdot\chi},
 \qquad
 T_gT_h=D_{\eta(g,h)}T_{gh}.
\end{equation}
Indeed, on a vector of degree \(a\), the first equality compares
\(\chi(\rho_g^{-1}a)\) on both sides. In the second equality,
the coherence scalar is evaluated in output degree
\(\rho_{gh}a\), so it acts after \(T_{gh}\).
Equations~\eqref{asf:eq-representation-operators} show that the map
\begin{equation}\label{asf:eq-group-representation}
 (\chi,g)\longmapsto D_\chi T_g
\end{equation}
is a representation of the group with multiplication
\eqref{asf:eq-extension-multiplication}.

Conversely, restrict a representation of \(K\) to its finite
abelian subgroup \(\widehat A\). It is a direct sum of character
eigenspaces. Through \(A\cong\widehat{\widehat A}\), this gives
the grading
\begin{equation}\label{asf:eq-weight-decomposition}
 V=\bigoplus_{a\in A}V_a,\qquad
 V_a=\{v:\chi v=\chi(a)v\text{ for every }\chi\in\widehat A\}.
\end{equation}
These eigenspaces are the \emph{weight spaces} for \(\widehat A\).
Let \(T_g\) be the operator of \((1,g)\).
The multiplication in \(K\) gives
\eqref{asf:eq-representation-operators}, so \(T_g\) regrades
according to \(\rho_g\) and defines \(u_g:F_g(V)\to V\)
satisfying~\eqref{asf:eq-equivariance}.
A linear map is a \(K\)-module map precisely when it preserves
these weight spaces and commutes with all \(T_g\).
Thus the two constructions are inverse also on morphisms.

Finally, in the strict model the equivariant tensor product has
operators \(T_g^V\otimes T_g^W\) on \(V\otimes W\), and its grading is obtained by adding
degrees. Consequently its operators \(D_\chi\) are
\(D_\chi^V\otimes D_\chi^W\). These are exactly the tensor product
operators for group representations. The unit and associativity
constraints agree under the same identification. We have therefore
constructed the tensor equivalence~\eqref{asf:eq-reconstructed-category}.
If the two group orders are coprime, the last assertion follows
from Lemma~\ref{asf:coprime-cohomology}.
\end{proof}
\begin{remark}\label{asf:remark-rebraiding}
Let \((A,q)\) be a pre-metric group of odd order, and
let a finite group \(G\) act on \(\C(A,q)\) by braided tensor
autoequivalences. The tensor equivalence of
Theorem~\ref{asf:odd-pointed-equivariantization} supplies a
Tannakian braiding on the same underlying equivariantization.
Explicitly, after transporting the action to the strict model,
replace \(c_q\) by the ordinary flip on \(\Vec_A\).
All regrading functors preserve this flip, and their coherence
maps remain monoidal. The flip therefore induces a symmetric
braiding on the equivariantization, corresponding to the usual
one on \(\Rep(K)\). The original braiding still records the quadratic form.
\end{remark}

\subsection{The order-two factor}

\begin{lemma}
\label{asf:two-primary-splitting}
Let \((A,q)\) be a pre-metric group such that
\(|A|\) is even and \(4\nmid|A|\), and let a finite group \(G\)
of odd order act on \(\mathcal P=\C(A,q)\) by braided tensor
autoequivalences. Write \(A_2\) for the two-primary subgroup
and \(A_{\mathrm{odd}}\) for the subgroup of odd-order elements.
Then \(A_2\cong C_2\), the action restricts to
\(\C(A_{\mathrm{odd}},q|_{A_{\mathrm{odd}}})\), and
\begin{equation}\label{asf:eq-two-primary-splitting}
 \mathcal P^G\simeq_{\br}
 \C(A_2,q|_{A_2})\boxtimes
 \C(A_{\mathrm{odd}},q|_{A_{\mathrm{odd}}})^G.
\end{equation}
The second factor is tensor equivalent to \(\Rep(H)\) for a
finite group \(H\) of order \(|A_{\mathrm{odd}}||G|\).
\end{lemma}
\begin{proof}
The primary decomposition is orthogonal for \(b_q\) by the coprimality of the orders of distinct primary summands. It gives a braided
equivalence
\(\mathcal P\simeq_{\br}
\C(A_2,q|_{A_2})\boxtimes
\C(A_{\mathrm{odd}},q|_{A_{\mathrm{odd}}})\).
Every isometry preserves the primary subgroups, and its restriction
to \(A_2\cong C_2\) is the identity. In particular, the odd
subcategory is stable under the action, including the coherence maps.

Choose a fixed model on \(A_2\) and the invariant model of
Lemma~\ref{asf:odd-canonical-model} on \(A_{\mathrm{odd}}\).
An isometry is represented by the identity functor on the first
factor and strict regrading on the second. Every braided
autoequivalence inducing the same isometry is monoidally isomorphic
to this representative. To verify this directly, compare their
tensor constraints after identifying their underlying linear
functors. The ratio satisfies the ordinary \(2\)-cocycle identity,
since the associators cancel. The braided identities make the
ratio symmetric. Apply Lemma~\ref{asf:symmetric-cocycle}.

We may therefore transport the action to these representatives.
The remaining coherence maps are characters of
\(A_2\oplus A_{\mathrm{odd}}\), and decompose into characters
of the two factors. The \(A_2\)-components form a cocycle in
\(Z^2(G,\widehat A_2)\), with trivial action on
\(\widehat A_2\cong C_2\). Since \(|G|\) is odd, the coprime
coherence lemma (Lemma~\ref{asf:coprime-cohomology}) makes this
component a coboundary. Changing the action identifications
removes it.

The resulting action is the trivial categorical action on the
first factor and the induced categorical action on the second.
Equivariant objects consequently decompose by their first-factor
degree, and equivariance is entirely in the second factor.
This gives the braided equivalence
\eqref{asf:eq-two-primary-splitting}, with the componentwise tensor
constraints and braiding. The final assertion follows from
Theorem~\ref{asf:odd-pointed-equivariantization}.
\end{proof}

\begin{corollary}
\label{asf:pointed-quotient-parity}
Let \(\C\) be a braided fusion category with
\(\FPdim(\C)=2m\), where \(m\) is an odd positive integer.
Let \(\Rep(G)\) be the maximal Tannakian subcategory of
\(\mathcal Z_2(\C)\), and suppose that \(\C_G\) is pointed.
Then exactly one of the following alternatives holds:
\begin{enumerate}
\item \(\C\simeq_\otimes\Rep(K)\) for a finite group \(K\) with
      \(|K|=2m\).
\item \(\C\simeq_{\br}\Sem_\varepsilon\boxtimes\D\), where
      \(\varepsilon\in\{+,-\}\), \(\D\) is a braided fusion
      category, and \(\D\simeq_\otimes\Rep(H)\) for a finite
      group \(H\) with \(|H|=m\).
\end{enumerate}
Furthermore,
\begin{equation}\label{asf:eq-pointed-quotient-parity}
 \C\text{ admits a fiber functor}
 \quad\Longleftrightarrow\quad
 2\mid\FPdim\bigl(\mathcal Z_2(\C)\bigr).
\end{equation}
\end{corollary}
\begin{proof}
Write \(\C_G\simeq_{\br}\C(A,q)\). By
\eqref{asf:eq-central-equivariantization} and \eqref{asf:eq-center-quotient},
\[
 |A||G|=2m,\qquad
 \FPdim(\mathcal Z_2(\C))=|G|\FPdim(\mathcal Z_2(\C_G)).
\]
The quotient center is \(\Vec\) or \(\sVec\) by maximality of
\(\E_\C\). If \(|A|\) is odd, then \(|G|\) is even, so the center of
\(\C\) has even dimension. Theorem~\ref{asf:odd-pointed-equivariantization}
gives the first alternative with \(|K|=|A||G|=2m\).

Suppose now that \(|A|\) is even, so \(|G|\) is odd. Let \(s\)
generate \(A_2\cong C_2\), and let \(X_s\) be the simple object of
degree \(s\). It centralizes the odd-order factor by orthogonality.
Since \(b_q(s,s)=q(s)^{-2}\), the object \(X_s\) is transparent
exactly when \(q(s)=\pm1\).

If \(\mathcal Z_2(\C_G)\simeq_{\br}\sVec\), its fermion is \(X_s\),
the unique invertible simple of order two. In this case \(q(s)=-1\).
The \(A_2\)-factor is therefore \(\sVec\), whose underlying tensor
category is tensor equivalent to \(\Rep(C_2)\). The splitting
\eqref{asf:eq-two-primary-splitting} now gives
\(\C\simeq_\otimes\Rep(C_2\times H)\), where
\(|H|=|A_{\mathrm{odd}}||G|=m\), and
\(\FPdim(\mathcal Z_2(\C))=2|G|\) is even.
If \(\mathcal Z_2(\C_G)\simeq_{\br}\Vec\), then \(X_s\) is not
transparent. Since \(q(s)^4=q(2s)=1\)
by~\eqref{asf:eq-quadratic-identities}, we have \(q(s)=i\) or \(-i\).
The \(A_2\)-factor is therefore semion, and the same splitting gives
the second alternative with \(|H|=m\). In this case
\(\FPdim(\mathcal Z_2(\C))=|G|\) is odd.

A group representation category has its usual forgetful fiber
functor. In the second alternative, a fiber functor on \(\C\)
would restrict to one on its semion factor, contradicting
Lemma~\ref{asf:semion-associator}. This proves both the
parity criterion~\eqref{asf:eq-pointed-quotient-parity} and
the disjointness of the two alternatives.
\end{proof}
The semion in Corollary~\ref{asf:pointed-quotient-parity} can be
detected without assuming that the central de-equivariantization is
pointed. For modular categories the splitting is
\cite[Corollary~3.9]{CCP}.
\begin{theorem}
\label{asf:general-semion-splitting}
Let \(\C\) be a braided fusion category over an algebraically closed
field of characteristic zero with \(\FPdim(\C)=2m\), where \(m\) is
an odd positive integer. Then \(\FPdim(\mathcal Z_2(\C))\) is odd if
and only if there are a semion category \(\Sem_\varepsilon\) and a
braided fusion category \(\D\) of dimension \(m\) such that
\begin{equation}\label{asf:eq-general-semion-splitting}
 \C\simeq_{\br}\Sem_\varepsilon\boxtimes\D .
\end{equation}
In this case \(\C\) admits no fiber functor.
\end{theorem}
\begin{proof}
By Proposition~\ref{asf:integrality-input}, \(\C\) is integral.
Suppose that its M\"uger center has odd dimension. It is Tannakian
by \cite[Corollary~2.50(i)]{DGNO}, so write it as \(\Rep(G)\).
By~\eqref{asf:eq-subcategory-divisibility}, the odd integer \(|G|\)
divides \(m\). The central quotient \(\B=\C_G\) is integral and
nondegenerate by Section~\ref{asf:subsec-central}. Its dimension is
\(2(m/|G|)\), twice an odd integer, and the positive spherical
structure makes it modular.

By \cite[\S2 and Corollary~3.9]{CCP}, \(\B\) splits as a braided
product of a semion and an
odd-dimensional modular category. The second factor has no invertible
object of order two: such an object would generate a subcategory of
dimension two, contrary to~\eqref{asf:eq-subcategory-divisibility}.
Consequently, \(\B\) has a unique invertible simple \(s\) of order two.
It generates the semion factor \(\mathcal S\), and every braided
autoequivalence of \(\B\) preserves this subcategory.

Restrict the \(G\)-action to \(\mathcal S\). Apply
Lemma~\ref{asf:two-primary-splitting} to this stable pointed
subcategory, whose group of simple objects is \(C_2\).
Its odd part is \(\Vec\), so
\(\mathcal S^G\simeq_{\br}\mathcal S\boxtimes\Rep(G)\).
The full inclusion \(\mathcal S\subseteq\B\) induces a fully
faithful braided functor \(\mathcal S^G\to\B^G\simeq_{\br}\C\).
The image of the first factor is therefore a semion subcategory of
\(\C\). Denote this image by \(\mathcal T\). Since \(\mathcal T\) is
nondegenerate, \cite[Theorem~3.13]{DGNO} gives
\(\C\simeq_{\br}\mathcal T\boxtimes\mathcal T'_\C\).
In~\eqref{asf:eq-centralizer-dimension}, the intersection
\(\mathcal T\cap\mathcal Z_2(\C)\) is \(\Vec\), so the complementary
factor \(\D=\mathcal T'_\C\) has dimension \(2m/2=m\).

Conversely, in~\eqref{asf:eq-general-semion-splitting} the semion is
nondegenerate, so \(\mathcal Z_2(\C)\simeq_{\br}\mathcal Z_2(\D)\).
The category \(\D\) is integral by Proposition~\ref{asf:integrality-input}.
Its center has dimension dividing the odd integer \(m\) by
\eqref{asf:eq-subcategory-divisibility}. Finally, any fiber functor
on \(\C\) would restrict to one on its semion subcategory, contrary
to Lemma~\ref{asf:semion-associator}.
\end{proof}
\section{Tensor structure in almost square-free dimension}
\label{asf:sec-almost-squarefree}

Yu proved that integral braided fusion categories of dimension \(p^n d\)
are group-theoretical for every prime \(p\) and \(0\le n\le4\), when
\(d\) is square-free and coprime to \(p\). His proof also identifies the
pointed quotient by the maximal Tannakian subcategory of the M\"uger center
\cite[Corollary~3.13 and its proof]{Yu}. We give a proof of this
reduction by prime-power decomposition and use it to reconstruct the
underlying tensor category for odd \(p\).

\begin{proposition}
\label{asf:central-pointed-reduction}
Let \(p\) be a prime and let \(d\) be a square-free positive integer
coprime to \(p\). Let \(\C\) be an integral braided fusion category
with \(\FPdim(\C)=p^n d\), where \(0\le n\le4\).
Let \(\E_\C=\Rep(G)\) be the maximal Tannakian subcategory of
\(\mathcal Z_2(\C)\). Then the de-equivariantization
\(\B=\C_G\) is pointed, its M\"uger center is braided equivalent to
\(\Vec\) or \(\sVec\), and
\begin{equation}\label{asf:eq-asf-central-reduction}
 \C\simeq_{\br}\B^G,\qquad |G|\FPdim(\B)=p^n d.
\end{equation}
Consequently, \(\C\) is group-theoretical.
\end{proposition}
\begin{proof}
By Section~\ref{asf:subsec-central}, the quotient \(\B\) is integral,
satisfies~\eqref{asf:eq-asf-central-reduction}, and has center \(\Vec\)
or \(\sVec\). Its dimension divides \(p^n d\), so it is again of the
form \(p^a d'\), with \(0\le a\le n\) and \(d'\mid d\).
Ostrik--Yu's criterion \cite[Corollary~4.9]{OY} therefore gives
nilpotence of \(\B\).
The prime-power decomposition of a braided nilpotent fusion category
\cite[Theorem~1.1]{DGNOnil} gives
\begin{equation}\label{asf:eq-prime-power-decomposition}
 \B\simeq_{\br}\boxtimes_{\ell\mid\FPdim(\B)}\B_\ell,
 \qquad \FPdim(\B_\ell)=\ell^{v_\ell(\FPdim(\B))},
\end{equation}
where \(\ell\) ranges over the prime divisors and \(v_\ell\) is the
\(\ell\)-adic valuation. Each factor is a fusion subcategory of
\(\B\), hence integral.
The M\"uger center of the product is the product of the M\"uger
centers of its factors. Thus every odd-primary factor is
nondegenerate, and a two-primary factor, if present, is nondegenerate
or slightly degenerate.

Equip each nondegenerate factor \(\B_\ell\) with its positive spherical
structure and pass to the complex model of
Section~\ref{asf:subsec-dimensions}. Write \(\FPdim(\B_\ell)=\ell^a\),
where \(1\le a\le4\). For odd \(\ell\), nonpointedness would imply
\(\ell^2\le\FPdim((\B_\ell)_{\mathrm{pt}})\le\ell^{a-3}\)
by \cite[Lemma~3.2]{CPS}, forcing \(a\ge5\). For \(\ell=2\), every
nonpointed possibility in this dimension range in
\cite[Proposition~3.1]{CPS} contains an Ising factor, hence a simple
of dimension \(\sqrt2\). This contradicts integrality. Thus each
nondegenerate primary factor is pointed, also over \(k\).

If the two-primary factor is slightly degenerate, its dimension is at
most \(2^4\). By \cite[Lemma~4.8]{Yu}, an integral nonpointed slightly
degenerate category of two-power dimension has dimension at least
\(2^6\). Hence this factor is pointed as well. The decomposition
\eqref{asf:eq-prime-power-decomposition} now shows that \(\B\) is
pointed. The criterion of Naidu--Nikshych--Witherspoon
\cite[Theorem~7.2]{NNW} gives group-theoreticality of \(\C\).
\end{proof}

\begin{remark}\label{asf:central-versus-crossed-reduction}
Maximality in Proposition~\ref{asf:central-pointed-reduction} is taken
inside \(\mathcal Z_2(\C)\). To see the distinction, take an odd prime
\(p\), put \(A=C_{p^2}\), choose a primitive \(p^2\)-th root \(\zeta\),
and set \(q(x)=\zeta^{x^2}\). The bicharacter
\(b_q(x,y)=\zeta^{2xy}\) is nondegenerate, so the M\"uger center of
\(\C(A,q)\) is \(\Vec\). On \(pA\), the form \(q\) is identically one.
Moreover, \(b_q(x,p)=1\) exactly when \(p\mid x\), so
\((pA)^\perp=pA\). Hence the pointed subcategory on \(pA\) is
Tannakian and is its own centralizer. It is therefore maximal
Tannakian in \(\C(A,q)\), although it has dimension \(p\) and the
maximal central Tannakian subcategory has dimension one.

The nilpotence criterion used in
Proposition~\ref{asf:central-pointed-reduction} applies for every
\(n\). The bound \(n\le4\) enters when proving pointedness of the
primary factors.
\end{remark}

\subsection{Tensor reconstruction for an odd prime}
\label{asf:subsec-main-proof}

\begin{proof}[Proof of Theorem~\ref{asf:odd-asf-tensor}]
Since \(4\nmid p^n d\), Proposition~\ref{asf:integrality-input}
gives integrality. By the pointed central reduction
(Proposition~\ref{asf:central-pointed-reduction}), the maximal
Tannakian subcategory \(\Rep(G)\) of \(\mathcal Z_2(\C)\) gives
\(\C\simeq_{\br}\C(A,q)^G\), where \(A\) is the group of simple
degrees in \(\C_G\) and \(|A||G|=p^n d\).
If \(p^n d\) is odd, then \(|A|\) is odd, so the reconstruction
theorem for odd pointed equivariantizations
(Theorem~\ref{asf:odd-pointed-equivariantization}) gives the first
assertion, including the order of \(K\). If \(p^n d\) is even,
it is twice an odd integer. The parity criterion for pointed central
de-equivariantizations (Corollary~\ref{asf:pointed-quotient-parity}) gives the two
disjoint alternatives and~\eqref{asf:odd-asf-fiber-criterion}.
In the semion case, Lemma~\ref{asf:semion-associator} gives the
underlying tensor equivalence~\eqref{asf:eq-asf-semion-tensor}.
\end{proof}

\begin{remark}
\label{asf:odd-asf-prior-work}
In odd dimension \(p^n d\) with \(n\le3\), Dong and Dai describe
semisimple quasitriangular Hopf algebra realizations by function
algebras or crossed products \cite{DD}. The reconstruction above
compares either realization with a group algebra. Write
\(\E_\C=\Rep(G)\) and \(\C_G\simeq_{\br}\C(A,q)\). Then
Theorem~\ref{asf:odd-pointed-equivariantization} gives an extension
\(1\to\widehat A\to K\to G\to1\). A normalized section \(s:G\to K\)
has factor set \(\eta(g,h)=s(g)s(h)s(gh)^{-1}\). Under the Fourier
isomorphism \(k\widehat A\cong k^A\), where \(k^A\) is the function
algebra on \(A\), formula~\eqref{asf:eq-extension-multiplication}
presents \(kK\) as a crossed product of \(k^A\) by \(kG\).
In particular, the factor set can give a nonsplit group extension,
as Example~\ref{asf:nontrivial-coherence-example} illustrates.

If a semisimple Hopf algebra \(\mathsf H\) realizes \(\C\), transport
its forgetful fiber functor along \(\C\simeq_\otimes\Rep(K)\).
By \cite[Proposition~5.14.4]{EGNO}, the transported functor corresponds
to a Drinfeld twist of \(kK\) up to gauge equivalence. Hopf
reconstruction therefore identifies \(\mathsf H\) with a twist of \(kK\).
When \(\C\) is pointed, \(K\) is abelian, and \(kK\cong k^{\widehat K}\).
\end{remark}

\subsection{Sharpness at exponent five}
\label{asf:sharpness-exponent-five}
We work over \(\mathbb C\) in this subsection. The following theorem shows that the bound \(n\le 4\) in
Theorem~\ref{asf:odd-asf-tensor} is sharp already for \(d=1\).

\begin{theorem}
\label{asf:sharpness-theorem}
For every odd prime \(p\) there is an integral
modular category \(\mathcal M_p\) of dimension \(p^5\) whose underlying
tensor category is inequivalent to \(\Rep(K)\) for every finite group
\(K\).
\end{theorem}
\begin{proof}
Let \(A=\mathbb F_p^3\), fix a primitive \(p\)-th root of unity \(\zeta\),
and put \(q(x,y,z)=\zeta^{y^2-4xz}\). The bilinear form in the exponent
of the polarization \(b_q\) has matrix
\[
 \begin{pmatrix}0&0&-4\\0&2&0\\-4&0&0\end{pmatrix}.
\]
Its determinant is \(-32\ne0\) in \(\mathbb F_p\), so
\(\mathcal P=\C(A,q)\) is pointed and modular. Identify \(G=C_p\) with
\((\mathbb F_p,+)\) and set
\[
 T_t(x,y,z)=(x,\ y+2tx,\ z+ty+t^2x)\qquad(t\in G).
\]
Then \(T_sT_t=T_{s+t}\), and
\((y+2tx)^2-4x(z+ty+t^2x)=y^2-4xz\). By
Lemma~\ref{asf:odd-canonical-model}, the \(T_t\) act by strict braided
regrading on the invariant model of \(\mathcal P\). Equipping the
graded spaces with Hilbert space structures gives a unitary model:
the associator is trivial and the braiding scalars have modulus one.
The strict action provides a categorical lift, so the first gauging
obstruction vanishes. For this lift the remaining obstruction lies in
\(H^4(G,U(1))\), with the trivial action on the group \(U(1)\) of
unit complex numbers. The cyclic cohomology formula
\cite[Theorem~6.2.2]{W} gives
\[
 H^4(C_p,U(1))\cong U(1)/(U(1))^p=0,
\]
since every unit complex number has a \(p\)th root in \(U(1)\).
Thus the lift extends to a braided \(G\)-crossed fusion category
\(\D=\bigoplus_{t\in G}\D_t\) with neutral component
\(\D_0=\mathcal P\) and all components nonzero
\cite[\S3]{CGPW}. Equivariantizing the entire extension yields a
modular category \(\mathcal M_p=\D^G\) of dimension
\(|G|^2\FPdim(\mathcal P)=p^5\) by \cite[Theorem~2]{CGPW}.
It is integral by Proposition~\ref{asf:integrality-input}.

To see a noninvertible simple explicitly, write \(X_a\) for the simple
of \(\mathcal P\) of degree \(a\in A\). The orbit of \((1,0,0)\) consists
of the \(p\) distinct points \((1,2t,t^2)\). Hence
\(\bigoplus_{t\in G}X_{(1,2t,t^2)}\), with \(G\) permuting the summands,
is an equivariant object in \(\mathcal P^G\). Its endomorphisms are
diagonal because the degrees are distinct, and equivariance makes the
diagonal entries constant on the orbit. Its endomorphism algebra is
therefore \(\mathbb C\). The full inclusion \(\mathcal P^G\subset\D^G\)
preserves this endomorphism algebra, giving a simple of dimension \(p\)
in \(\mathcal M_p\).

Suppose a tensor equivalence \(\mathcal M_p\simeq_\otimes\Rep(K)\)
existed. It would transport the modular braiding to \(\Rep(K)\).
By \cite[Example~2.1]{Nat06}, the corresponding \(R\)-matrix is
supported on \(\mathbb C\Gamma\otimes\mathbb C\Gamma\) for an abelian
normal subgroup \(\Gamma\) of \(K\). On every representation inflated
from \(K/\Gamma\), the group algebra \(\mathbb C\Gamma\) acts through
the counit. The two counit identities for \(R\) therefore make both
crossings with such a representation the ordinary flip. Thus
\(\Rep(K/\Gamma)\) lies in the M\"uger center. Nondegeneracy forces
\(K=\Gamma\), so every simple of \(\Rep(K)\) has dimension one,
contradicting the simple of dimension \(p\) above.
\end{proof}

By modularity the M\"uger center of \(\mathcal M_p\) is \(\Vec\), so its
maximal central de-equivariantization is \(\mathcal M_p\) itself. Thus
this family realizes the nonpointed case of
Proposition~\ref{asf:exponent-five-obstruction}\textup{(i)}.
At \(p=3\) it has nine invertible simples and twenty-six simples of
dimension three, the nonpointed type in \cite[Lemma~4.3]{XD}.

The same family shows sharpness for both group-representation conclusions of
Theorem~\ref{asf:odd-asf-tensor}\textup{(ii)} at \(d=2\).
The product \(\mathcal M_p\boxtimes\Rep(C_2)\), with the symmetric
braiding on \(\Rep(C_2)\), has M\"uger center of dimension two.
A tensor equivalence with \(\Rep(K)\) would identify its subcategory
\(\mathcal M_p\boxtimes\Vec\) with \(\Rep(K/N)\) for some normal
subgroup \(N\) of \(K\), since every fusion subcategory of \(\Rep(K)\)
has this form (take \(N\) to be the intersection of the kernels of
its objects). This contradicts Theorem~\ref{asf:sharpness-theorem}.
The product \(\Sem_+\boxtimes\mathcal M_p\) has trivial M\"uger center.
Since \(\mathcal M_p\) has odd dimension, it has no invertible object
of order two. The product therefore has a unique semion subcategory,
whose centralizer is \(\Vec\boxtimes\mathcal M_p\). Its odd-dimensional
complement consequently fails the group-representation conclusion
in Theorem~\ref{asf:odd-asf-tensor}\textup{(ii)(b)}.

\subsection{Coherence and the prime two}

\begin{example}
\label{asf:nontrivial-coherence-example}
Let \(p\) be an odd prime, and let \(\zeta\in k\) be a
primitive \(p\)th root of unity.
On \(A=C_{p^2}\), define \(q(x)=\zeta^{x^2}\).
Then
\begin{equation}\label{asf:eq-degenerate-cyclic-form}
 b_q(x,y)=\zeta^{2xy},\qquad
 A^\perp=pC_{p^2},\qquad q|_{pC_{p^2}}=1.
\end{equation}
The M\"uger center of \(\C(A,q)\) is the Tannakian subcategory
on \(pA\), of dimension \(p\). De-equivariantizing \(\C(A,q)\)
by this center gives \(\C(A/pA,\bar q)\), where \(A/pA\cong C_p\) by reduction
modulo \(p\) and \(\bar q(x+pA)=\zeta^{x^2}\).
To see this explicitly, choose the trivial-associator model of
Lemma~\ref{asf:odd-canonical-model}. The regular algebra is the sum of
the lines in degrees \(pA\), with multiplication adding these degrees
after a choice of character basis. Multiplication by a basis element
in degree \(h\in pA\) identifies the module components in degrees
\(x\) and \(x+h\). Translation acts freely on each coset, so a simple
module has one-dimensional components on one coset and vanishes on
the others. It is the free module on a simple of degree \(x\).
The free-module functor is braided and therefore gives its coset
the self-braiding \(q(x)=\bar q(x+pA)\).

The induced \(C_p\)-action on the simple objects of the
de-equivariantization is trivial, since \(O(C_p,\bar q)=C_2\).
The categorical action is not trivial. If it were,
its equivariantization would be
\(\C(C_p,\bar q)\boxtimes\Rep(C_p)\), whose group of
invertible objects is \(C_p\times C_p\).
The original category has group of invertible objects
\(C_{p^2}\), whereas \(C_p\times C_p\) has exponent \(p\).

This coherence can be described explicitly.
Represent elements of \(G=C_p\) by \(0,\ldots,p-1\),
identify \(t\in C_p\) with the character \(a\mapsto\zeta^{ta}\)
of \(C_p\), and take
\begin{equation}\label{asf:eq-carry-cocycle}
 \eta(g,h)(a)=\zeta^{\,a\left\lfloor(g+h)/p\right\rfloor}.
\end{equation}
The floor is the integer part. The identity expressing
equality of the total carry in adding three integers is
the \(2\)-cocycle identity for~\eqref{asf:eq-carry-cocycle}.
In the extension multiplication
\eqref{asf:eq-extension-multiplication}, the element
corresponding to \((0,1)\in C_p\times C_p\) has
\(p\)th power \((1,0)\), and has order \(p^2\).
Thus the extension is
\begin{equation}\label{asf:eq-cyclic-nonsplit-extension}
 1\longrightarrow C_p\longrightarrow C_{p^2}
 \longrightarrow C_p\longrightarrow1.
\end{equation}
Choose a primitive \(p^2\)th root \(\xi\) with \(\xi^p=\zeta\).
The character indexed by \(x\in C_{p^2}\) sends the generator
\((0,1)\) to \(\xi^x\) and its \(p\)th power to \(\zeta^x\).
Its restriction to the kernel has weight \(x\bmod p\), so
the induced self-braiding is \(\zeta^{x^2}\). Thus this choice of
cocycle reconstructs the given pre-metric group and realizes a
nonzero extension class.
\end{example}

\begin{example}
\label{asf:cyclic-four-example}
Let \(\xi\in k\) be a primitive eighth root of unity, and let
\(\B=\C(C_4,q)\), where \(q(x)=\xi^{x^2}\).
The form is well-defined modulo four because
\((x+4)^2-x^2\) is divisible by eight.
Its bicharacter is \(b_q(x,y)=\xi^{2xy}\), which has
trivial radical. Hence \(\B\) is an integral pointed
nondegenerate braided category of dimension four.
The unique nonzero element of order two has
\(q(2)=\xi^4=-1\); in particular \(\B\) contains no
semion subcategory.

The underlying fusion category of \(\B\) has no fiber functor.
By \cite[Example~5.1.3]{EGNO}, a fiber functor on a pointed
fusion category forces its associator class to be trivial. Transporting
the braiding to the resulting \(\Vec_{C_4}\), whose associator is
trivial, would make it a bicharacter \(c\).
For a generator \(a\) of \(C_4\), such a bicharacter
satisfies \(c(a,a)^4=c(4a,a)=1\).
Here every generator has \(q(a)^4=-1\), a contradiction.

Thus \(\B\) is not tensor equivalent to \(\Rep(K)\).
Nor is it tensor equivalent to
\(\Vec_{C_2}^{\omega_{\mathrm{nt}}}\boxtimes\Rep(H)\):
dimension would force \(|H|=2\), and the latter category
has group of invertible objects \(C_2\times C_2\), whereas
\(\B\) has group \(C_4\).
Proposition~\ref{asf:central-pointed-reduction} still
applies, with trivial equivariantizing group.
\end{example}

Thus at the prime two the restriction to the order-two subgroup alone
does not determine the underlying tensor structure. For comparison,
give \(\Vec_{C_2\times C_2}\) the braiding
\(c((a,b),(u,v))=(-1)^{av}\). Its double braiding is
\((-1)^{av+bu}\), which is nondegenerate, and the element \((1,1)\)
has self-braiding \(-1\). This gives a second nondegenerate category
of dimension four with the same restricted form on a subgroup of order
two. Its trivial associator gives a fiber functor by forgetting the grading.
Conversely, different quadratic forms can give the same underlying tensor
category. Reversing the braiding in
Example~\ref{asf:cyclic-four-example} replaces \(q\) by \(q^{-1}\)
and leaves its associator fixed. The two forms are
nonisometric, since both generators of \(C_4\) have square \(1\)
modulo eight, and \(\xi\ne\xi^{-1}\).

\subsection{Applications to braidings and Hopf algebras}

\begin{corollary}
\label{asf:jordan-larson-no-braiding}
Let \(p,q\) be distinct primes, and let \(\C\) be an
integral fusion category with
\(\FPdim(\C)=p^2q\).
If \(\C\) is not group-theoretical, then it admits no braiding.
In particular this applies to the non-group-theoretical
categories in \cite[Theorem~1.1]{JL}.
If both \(p\) and \(q\) are odd, then \(\C\) admits a
braiding if and only if it is tensor equivalent to
\(\Rep(K)\) for some finite group \(K\) of order \(p^2q\).
\end{corollary}
\begin{proof}
Jordan and Larson denote the repeated prime by \(q\). In our notation,
a braiding on \(\C\) allows us to apply
Proposition~\ref{asf:central-pointed-reduction} with \(n=2\) and
\(d=q\). This gives group-theoreticality, as also follows directly
from \cite[Corollary~3.13]{Yu}. If both primes are odd,
Theorem~\ref{asf:odd-asf-tensor} gives a tensor equivalence with
\(\Rep(K)\) of the stated order. Conversely, such an equivalence
transports the usual symmetric braiding on \(\Rep(K)\) to \(\C\).
\end{proof}
\begin{remark}\label{byz}
Let \(\mathsf H\) be a semisimple quasitriangular Hopf algebra of
square-free dimension over \(k\).
By Corollary~\ref{asf:squarefree-tensor}, the forgetful fiber functor
on \(\Rep(\mathsf H)\) gives the group-representation alternative.
Transport the functor along a tensor equivalence
\(\Rep(\mathsf H)\simeq_\otimes\Rep(K)\). Fiber functors on \(\Rep(K)\)
correspond to conjugacy classes of subgroups \(L\le K\) equipped with
a nondegenerate \(2\)-cocycle class
\cite[Proposition~5.14.4 and Corollary~7.12.24]{EGNO}.
Here nondegeneracy means that the corresponding twisted group algebra
of \(L\) is a full matrix algebra over \(k\). Its dimension \(|L|\)
is therefore a square dividing the square-free integer \(|K|=\dim_k\mathsf H\).
Thus \(L=1\), the fiber functor is tensor isomorphic to the usual one,
and Hopf reconstruction gives \(\mathsf H\cong kK\).
\end{remark}
\subsection{A necessary condition in dimension \texorpdfstring{\(p^5d\)}{p to the fifth times d}}
\begin{proposition}
\label{asf:exponent-five-obstruction}
Let \(p\) be a prime and \(d\) a square-free positive integer coprime
to \(p\). Let \(\C\) be an integral braided fusion category over an
algebraically closed field of characteristic zero with
\(\FPdim(\C)=p^5d\). Let \(\Rep(G)\) be the maximal Tannakian
subcategory of its M\"uger center \(\mathcal Z_2(\C)\), and let
\(\B=\C_G\). If \(\B\) is nonpointed, then \(p\nmid |G|\) and
\(\B\) has a unique nonpointed primary factor \(\B_p\) in its
braided prime-power decomposition. It satisfies
\begin{equation}\label{asf:eq-exponent-five-primary}
 \FPdim(\B_p)=p^5,\qquad
 \B\simeq_{\br}\B_p\boxtimes\mathcal P,\qquad
 \FPdim(\mathcal P)=d/|G|,
\end{equation}
where \(\mathcal P\) is pointed.
\begin{enumerate}
\item If \(p\) is odd, then \(p\nmid\FPdim(\mathcal Z_2(\C))\)
and \(\B_p\) is modular for its positive spherical structure, with
\begin{equation}\label{asf:eq-exponent-five-modular}
 \operatorname{FSexp}(\B_p)\in\{p,p^2\},\qquad
 \FPdim((\B_p)_{\mathrm{pt}})=p^2 .
\end{equation}
The factor \(\B_p\) has exactly \(p^3-1\) noninvertible simple objects, each of
Frobenius--Perron dimension \(p\).
\item If \(p=2\), then \(\B_2\) is an integral nondegenerate braided
category of dimension \(32\), hence modular for its positive spherical
structure. The full quotient \(\B\) is nondegenerate and
\(\FPdim(\mathcal Z_2(\C))=|G|\) is odd.
\end{enumerate}
\end{proposition}
\begin{proof}
The category \(\B\) is integral and has M\"uger center \(\Vec\) or
\(\sVec\) by~\eqref{asf:eq-center-quotient} and maximality of
\(\Rep(G)\). Since \(|G|\FPdim(\B)=p^5d\), the nilpotence criterion
\cite[Corollary~4.9]{OY} applies to \(\B\).
If \(p\mid |G|\), then the \(p\)-adic valuation of \(\FPdim(\B)\)
is at most four and all other prime valuations are at most one.
Proposition~\ref{asf:central-pointed-reduction}, applied to \(\B\)
itself, gives pointedness because the maximal Tannakian subcategory of
\(\mathcal Z_2(\B)\) is trivial. Thus nonpointedness forces
\(p\nmid |G|\), and consequently \(|G|\mid d\).

In~\eqref{asf:eq-prime-power-decomposition}, each factor at a prime
\(\ell\ne p\) has dimension \(\ell\) and has center \(\Vec\) or
\(\sVec\). Its maximal central Tannakian subcategory is therefore
trivial. Proposition~\ref{asf:central-pointed-reduction} makes this
factor pointed. Their product is the factor \(\mathcal P\) in
\eqref{asf:eq-exponent-five-primary}, so \(\B_p\) must be nonpointed.

For odd \(p\), the center of \(\B_p\) has dimension at most two
and divides \(p^5\), so it is trivial. Thus \(\B_p\) is modular for
its positive spherical structure. Also,
\(\FPdim(\mathcal Z_2(\C))=|G|\FPdim(\mathcal Z_2(\B))\) by
\eqref{asf:eq-central-equivariantization}--\eqref{asf:eq-center-quotient},
and both factors are prime to \(p\). Apply \cite[Lemma~3.2 and
Proposition~3.3]{CPS} to the complex model of \(\B_p\) from
Section~\ref{asf:subsec-dimensions}. The lemma gives
\(p^2\le\FPdim((\B_p)_{\mathrm{pt}})\le p^{5-3}=p^2\), and the
proposition gives \(\operatorname{FSexp}(\B_p)\in\{p,p^2\}\).
Both numerical conclusions hold over \(k\).

The universal grading of \(\B_p\) therefore has \(p^2\) components,
each of dimension \(p^3\). The dimension of each simple is a power of
\(p\), since it is an integer whose square divides \(p^5\)
\cite[Proposition~8.14.6]{EGNO}. A simple of dimension at least
\(p^2\) would contribute at least \(p^4\) to a single component.
Every noninvertible simple consequently has dimension \(p\), and
\eqref{asf:eq-category-dimension} gives \((p^5-p^2)/p^2=p^3-1\)
such simples.

For \(p=2\), the nonpointed integral factor \(\B_2\) has dimension
\(2^5=32\). Its center is either \(\Vec\) or \(\sVec\). By
\cite[Lemma~4.8]{Yu}, a nonpointed integral slightly degenerate
category of dimension \(2^a\) requires \(a\ge6\). Thus \(\B_2\) is
nondegenerate. Each odd-primary factor also has trivial center, since
its center has dimension at most two and divides its odd dimension.
Taking the product shows that \(\mathcal Z_2(\B)\simeq_{\br}\Vec\).
The center dimension formula above now gives
\(\FPdim(\mathcal Z_2(\C))=|G|\), which is odd because \(|G|\mid d\).
\end{proof}
In particular, let \(p\) be odd and let \(\C\) be a braided fusion
category of dimension \(p^5d\) over \(k\), with \(d\) square-free and
prime to \(p\). Put \(\Rep(G)=\E_\C\) and \(\B=\C_G\). Integrality
follows from Proposition~\ref{asf:integrality-input}.
If \(p\mid\FPdim(\mathcal Z_2(\C))\), then \(\B\) is pointed by
Proposition~\ref{asf:exponent-five-obstruction}. If the center
dimension is prime to \(p\), then \(|G|\) divides that dimension,
so \(p\nmid |G|\). Hence \(\FPdim(\B)=p^5d/|G|\) has
\(p\)-adic valuation five. As in the proof above, \(\B\) is
nilpotent, and its \(p\)-primary factor has dimension \(p^5\) and
trivial center. Thus \(\B_p\) is modular whether or not \(\B\) is
pointed. By part~\textup{(i)}, an exponent
\(\operatorname{FSexp}(\B_p)\notin\{p,p^2\}\) forces \(\B\) to be pointed.
Under either sufficient condition,
Theorem~\ref{asf:odd-pointed-equivariantization} in odd dimension and Corollary~\ref{asf:pointed-quotient-parity} in even dimension
give the tensor alternatives of Theorem~\ref{asf:odd-asf-tensor}.

\section*{Acknowledgments}
The first author was supported by a grant of the Ministry of Research,
Innovation and Digitization, CNCS--UEFISCDI, project number
PN-IV-P1-PCE-2023-0264, within PNCDI IV.

ChatGPT (OpenAI) was used to examine examples and to assist with
drafting and revising portions of the manuscript. The authors take
responsibility for the statements, proofs, references, and final text.

\end{document}